\documentclass[12pt, amstex, amssymb]{article}
\title{Classical lambda calculus in modern dress}
\author{J. M. E. Hyland}
\input{diagrams.sty}
\newtheorem{theorem}{Theorem}[section]
\newtheorem{proposition}[theorem]{Proposition}

\newtheorem{lemma}[theorem]{Lemma}
\newtheorem{corollary}[theorem]{Corollary}
\newtheorem{definition}[theorem]{Definition}
\newcommand{\lscott}{[\![}
\newcommand{\rscott}{]\!]}
\newcommand{\ret}{\triangleleft}
\newcommand{\mcat}{\cal}
\newcommand{\cat}{\mathbf}

\newcommand{\app}{{\bf 1}}

\newcommand{\cI}{{\bf I}}
\newarrow{Eq}{=}{=}{=}{=}{=}
\usepackage{amsmath}

\begin{document}
\maketitle
\begin{abstract}
Recent developments in the categorical foundations of
universal algebra have given an impetus to an 
understanding of the lambda calculus coming from 
categorical logic: an interpretation is a semi-closed 
algebraic theory. Scott's representation 
theorem is then completely natural and leads to a
precise Fundamental Theorem showing the essential equivalence 
between the categorical and more familiar notions. 
\end{abstract}
\section{Introduction}

The $\lambda$-calculus is one of the great discoveries of logic
in the 20th century, but the question of its semantics has proved
vexed. Barendregt's monumental text \cite{Bar81} offered a variety
of approaches and it is telling that under the influence of Scott
\cite{Scott80} the treatment was largely rewritten
for the revised edition \cite{Bar84}. I offer a completely
new approach to the semantics. I believe that it shows that
the $\lambda$-calculus is a simple natural object of 
mathematical study.

The definition explained here is that an 
interpretation
of the lambda calculus is an algebraic theory equipped
with semi-closed structure: I call such a structure
a $\lambda$-theory. In discussing semantics,
I use the neutral
term ``interpretation'' in order to respect as
far as possible established usage of the
terms ``$\lambda$-algebra''
and ``$\lambda$-model''. One benefit of the view which I propose 
is that the inductive definition of an interpretation is done 
once and for all in the abstract setting: there is no further need
for it in individual cases. That is characteristic of categorical logic.
An algebraic theory is admittedly a slightly 
more complicated mathematical structure than one expects
in semantics, but interpreting abstraction involves
handling free variables and the algebraic theory makes them
explicit. Moreover the definition can be used effectively
in giving examples: common interpretations are naturally
presented as semi-closed algebraic theories so there is 
a real gain. A further benefit of
making the notion of theory central is that it leads quickly 
to fundamental results.
Specifically I show that Scott's interpretation
by reflexive objects in cartesian closed categories
arises naturally.

The initial $\lambda$-theory $\Lambda$ has a presentation
by the syntax of the $\lambda$-calculus. As an algebraic theory
it has algebras, and
a $\Lambda$-algebra is a clean version of a
valuation interpretation or environment model 
of the $\lambda$-calculus, bypassing the issue of
weak extensionality. 
The initial $\Lambda$-algebra is the closed term
interpretation traditionally written $\Lambda_0$.
Remarkably every $\lambda$-theory
is the theory of extensions of some $\Lambda$-algebra.
This gives a very tight equivalence between
the categories of $\lambda$-theories and that
of $\Lambda$-algebras. I want to call this the Fundamental 
Theorem of the Lambda Calculus.

I started writing this paper with a section on Combinatory
Algebra and $\lambda$-calculus, justifying 
the $\lambda$-theory definition. Though there is a close link
with the Fundamental Theorem I have cut that material
for reasons of space. I hope to present an account soon.
Hints as to what is involved are in \cite{Freyd89}
and especially \cite{Sel02}, where what I
call the $\lambda$-theory of a $\Lambda$-algebra
is treated albeit from a completly different point of view.
There are other omissions. In particular I am sorry not
to make more links with the universal algebra aspects of
the $\lambda$-calculus pioneered by Antonino Salibra and his
co-workers, see \cite{MS06}.
 
The notion of interpretation of the lambda calculus which I 
propose derives from an approach to semantics pioneered 
long ago by Andrew Pitts: one systematically takes contexts seriously.
This had a major unrecognised influence on early thinking about
categorical models for Linear Logic \cite{BBPHa}
and the idea to consider theories with structure
dates to that era. However only more recently have 
I realised how well that works for the $\lambda$-calculus.
The multicategory theory perspective now seems compelling in the light of
the foundations for notions of algebra \cite{Hyl13} provided by
Kleisli Bicategories \cite{FGHW}. I refrain from
discussing that background and also say nothing 
about the related theory of variable binding initiated in
\cite{FPT99} and extended variously for example in
\cite{TP06a}. While this is important, 
it is unnecessary for a first
appreciation of the $\lambda$-calculus.

This paper was refereed in an intelligent and helpful way. I hope that
I have profited from advice to clarify the exposition and to
give more precise details. Even so my account is very abstract
and more condensed than I would have liked. I have 
written it to honour Corrado B\"{o}hm, the creator of the 
technique \cite{Bohm68} which most influenced my understanding 
of the subject.  I hope that readers will find it suggestive
for further research. In this foundational treatment I cannot 
get as far as B\"{o}hm's Theorem but I
close the paper with suggestions for the future
including reflections on B\"{o}hm's legacy.

\section{Algebraic theories}

\subsection{Algebraic theories as cartesian operads}
An algebraic theory is a theory of equality on
terms. A clean mathematical expression of the idea is 
via multicategories or operads, so
abstractly via monads in some Kleisli bicategory \cite{FGHW}. I
only need
the concrete version of a one object cartesian multicategory or 
cartesian operad. Write ${\bf Sets}$ for the category of sets
and ${\cat F}$ for a standard skeleton 
of finite sets.
\begin{definition}
An {\em algebraic theory} $\mcat T$ is first
a functor ${\mcat T}: {\cat F} \to {\bf Sets}$:
so we have sets ${\cal T}(n)$ 
of $n$-ary multimaps with variable renamings. In addition, $\mcat T$
is equipped with projections ${\bf pr}_1, \cdots ,{\bf pr}_n \in {\mcat T}(n)$
including as special case the identity ${\bf id} \in {\cal T}(1)$. Finally
there are compositions ${\cal T}(n) \times {\cal T}(m)^n  \to 
{\cal T}(m)$ which are associative, unital, compatible with projections
and natural in $n$ and $m$ (or dinatural in $m$).
A {\em map $F : {\mcat S} \to {\mcat T}$ of algebraic theories}
is a natural transformation with components
$F_n : {\mcat S}(n) \to {\mcat T}(n)$
preserving projections and composition.
\end{definition}
Clearly from the definition we get a category of algebraic theories.
It is easily seen to be locally finitely presentable but I shall not
need that fact.

I have defined an algebraic theory in modern style as
what could be called a cartesian operad. Experts will be 
aware of other formulations in the operads tradition, and 
I comment later on the relation with Lawvere theories and monads.
But in essence the notion of algebraic theory
is the same as that of abstract clone \cite{Tay93} familiar in universal
algebra. Abstract clones appear without the name 
already in \cite{Tay73}. The functorial perspective
is simply a clean way to handle variable reindexing.
Note in particular the elements of ${\mcat T}(0)$ which are the constants
of the theory. (These are often omitted in clone theory.)
We have unique maps $0 \to n$ so we can identify the constants within
each ${\mcat T}(n)$. I shall take that kind of thing for
granted hereafter.

I make a brief remark on the syntactic point of view.
Writing $\Gamma \vdash t$ for $t$ a term with variable
declaration $\Gamma$, the definition of algebraic theory
encapsulates the basic 
principles of term formation
\[
 \dfrac{}{\Gamma \vdash x}
 \hspace{1.5cm} 
\dfrac{\Gamma \vdash t \quad \Delta \vdash s_1 , \cdots 
 \Delta \vdash s_n }{\Delta \vdash t({\bf s})}
\]
with $x$ declared in $\Gamma$ and $t({\bf s})= t(s_1, \cdots s_n)$ the result
of substituting the string of terms $\bf s$ in $t$. An algebraic
theory or abstract clone can be presented by allowing equations
between terms. The basic
syntactic rule of equality, the substitution of equals for 
equals, is implicit. 
I imagine that readers will have no problem with algebraic
theories arising thus
from concrete syntax. However it is good to be aware of another
source of examples. Suppose that $\cal C$ is a category with 
products and $X$ an object of $\cal C$. Then 
using the evident composition
${\cal C}(X^n,X) \times {\cal C}(X^m,X)^n \to {\cal C}(X^m,X)$,
one has
an algebraic theory, the {\em endomorphism theory} of $X$,
 with underlying functor
${\cal C}(X^n,X)$.

\subsection{Algebras for theories}\label{algebras}

Very general notions of interpretation for an algebraic theory 
are supported by \cite{FGHW} but we only need
the basic interpretation in the category $\bf Sets$ of sets. 
\begin{definition}
An algebra for an algebraic theory
$\mcat T$ is a set $A$ with an associative unital action
${\mcat T}(n) \times A^n \to A$ of $\mcat T$,
natural in $n$. If $A$ and $B$ are ${\mcat T}$-algebras then a homomorphism
from $A$ to $B$ consists of a map $f: A \to B$ respecting the actions
in the sense that the evident diagram commutes.
\end{definition}
Concretely one can take unital associative to mean first that
the projections ${\bf pr}_i \in {\mcat T}(n)$ act as projections
and that the two evident maps 
${\mcat T}(n) \times {\mcat T}(m)^n \times A^m \to A$
are equal.

From the definition we get a category $\mcat T$-algebras which I
write ${\rm Alg}({\mcat T})$. 
I give some background on this category. 
Further details can be extracted from
\cite{ARV11}. First note that
the compositions
${\mcat T}(m) \times {\mcat T}(n)^m \to {\mcat T}(n)$
give each ${\mcat T}(n)$ the structure of a $\mcat T$-algebra.
\begin{proposition}\label{free}
The algebra ${\mcat T}(n)$
is the free $\mcat T$-algebra on $n$ generators; that is,
for $A \in {\rm Alg}({\mcat T})$, we have
${\rm Alg}({\mcat T})({\mcat T}(n), A) \cong A^n$, natural in
$A$ and $n$.
\end{proposition}
{\bf Proof} This is a form of the Yoneda Lemma: the coordinates
in $A^n$ come from the $n$ projections 
${\bf pr}_1, \cdots ,{\bf pr}_n \in {\mcat T}(n)$.

\begin{corollary}\label{corolcoprod} ${\mcat T}(0)$ is the initial algebra and
${\mcat T}(n) + {\mcat T}(m) \cong {\mcat T}(n+m)$ gives binary coproducts.
\end{corollary}
The main closure properties of ${\rm Alg}({\mcat T})$ follow as in
\cite{ARV11}.
\begin{proposition}
The category ${\rm Alg}({\mcat T})$ of $\mcat T$-algebras is complete
and cocomplete.
\end{proposition}
{\bf Proof}
Completeness is evident as the forgetful functor to $\bf Sets$
creates limits. It also creates sifted colimits so it suffices
to define finite coproducts. The corollary above gives
the initial algebra and also helps justify the description
of the coproduct $A + B$ as a sifted colimit
$A + B = \int^{m,n} {\rm Alg}({\mcat T})({\mcat T}(m), A)
\times {\rm Alg}({\mcat T})({\mcat T}(n), B) \times{\mcat T}(m+n)$.\\[0.4em]
It is clear that if $F: {\mcat S} \to {\mcat T}$
is a map of algebraic theories then composition gives a
functor $F^*:{\rm Alg}({\mcat T}) \to {\rm Alg}({\mcat S})$.
I shall not need its left adjoint, though there are hints of it
in what follows. For each $n$ there is a map of
$\mcat S$-algebras ${\mcat S}(n) \to F^*{\mcat T}(n)$ with underlying
map $F_n$ carrying the $n$ generators
to the $n$ generators.  Let
$B$ be a $\mcat T$-algebra. It is easy to see that composing the
action of $F^*$ with the map induced by ${\mcat S}(n) \to F^*{\mcat T}(n)$
in
\[
B^n \cong {\bf Alg}({\mcat T})({\mcat T}(n),B) \longrightarrow
{\bf Alg}({\mcat S})(F^*{\mcat T}(n),F^*B) \longrightarrow
{\bf Alg}({\mcat S})({\mcat S}(n),F^*B) \cong B^n
\]
gives the identity on the set $B^n$.
\begin{proposition}\label{equiv}
Suppose that the functor $F^* :
{\rm Alg}({\mcat T}) \to {\rm Alg}({\mcat S})$ is an equivalence of categories.
Then $F$ is an isomorphism of algebraic theories.
\end{proposition}
{\bf Proof}
The first arrow above is an isomorphism as $F^*$ is full and faithful.
$F^*$ is essentially surjective on objects so we can put
any $\mcat S$-algebra $A$ in place of $F^*B$. Hence 
${\bf Alg}({\mcat S})(F^*{\mcat T}(n),A) \cong
{\bf Alg}({\mcat S})({\mcat S}(n),A)$ and $F^*$
full and faithful gives naturality in $A$. Thus 
the ${\mcat S}(n) \to F^*{\mcat T}(n)$ are isomorphisms.
Hence so are the $F_n : {\mcat S}(n) \to {\mcat T}(n)$, that is,
$F$ is an isomorphism.\\[0.4em]  
Note that in this situation $F^*$ is in
fact necessarily an isomorphism of categories.

Returning now to coproducts, $A[n] = A + {\mcat T}(n)$ gives the free extension
of a $\mcat T$-algebra $A$ by $n$ indeterminates.
Let $A$ be a model of an algebraic theory $\mcat S$.
There is an algebraic theory
${\mcat S}_A$ whose algebras are $\mcat S$-algebras
equipped with an ${\mcat S}$-algebra map from $A$. 
As sets ${\mcat S}_A(n) = A[n] = A + {\mcat S}(n)$. Giving 
the structure of an algebraic theory is routine.
(In case $A = {\mcat S}(p)$, Corollary \ref{corolcoprod}
gives the simple description
${\mcat S}_{{\mcat S}(p)}(n) = {\mcat S}(n+p)$ with action
not affecting the parameters in $p$.)
We have an evident map ${\mcat S} \to {\mcat S}_A$ inducing a functor
${\rm Alg}({\mcat S}_A) \to {\rm Alg}({\mcat S})$ and 
an isomorphism between $A$ and ${\mcat S}_A(0)$ regarded
as an $\mcat S$-algebra. A concrete way to think about
${\mcat S}_A$ 
is that it is obtained from $\mcat S$
by adding constants for the new elements of $A$
and equations coming from $A$
but no other equations.

Let $F: {\mcat S} \to {\mcat T}$ be a map of algebraic theories.
Set $A = F^*{\mcat T}(0)$ and observe that $F$  
factors through ${\mcat S} \to {\mcat S}_A$ via a comparison
${\mcat S}_A \to {\mcat T}$. If the latter is
an isomorphism $\mcat T$ is {\em the theory
of extensions of a model of} $\mcat S$. 
Proposition \ref{equiv} implies the following.
\begin{proposition}\label{algext}
$F: {\mcat S} \to {\mcat T}$ is a theory of extensions of a model
if and only if ${\mcat S}_A \to {\mcat T}$ induces an equivalence
between ${\rm Alg}({\mcat T})$ and the coslice
category $A / {\rm Alg}({\mcat S})$.
\end{proposition}
Whenever $F : {\mcat S} \to {\mcat T}$ 
is a map of theories and $B$ a $\mcat T$-algebra,
we get a map ${\mcat S}_{F^*B} \to {\mcat T}_B$ of theories.
If $\mcat T$ is a theory of extensions of a model of $\mcat S$,
then this will be an isomorphism.

\subsection{The presheaf topos}\label{topos}

A $\mcat T$-algebra is a set $A$ equipped with an action
of $\mcat T$ on the left. But $\mcat T$ can also act on
the right. In the operad literature one talks 
of a module: I prefer to say presheaf.
\begin{definition}
A {\em presheaf} $X$ over an algebraic theory $\mcat T$ is a functor
$X: {\cat F} \to {\bf Sets}$ equipped with an action
$X(m) \times {\mcat T}(n)^m \to X(n)$
compatible with the operations of $\mcat T$. A {\em map of presheaves}
is a functor commuting with the action of $\mcat T$.
\end{definition}
We get a category of presheaves over $\mcat T$ which I write
$P{\mcat T}$. Evidently $\mcat T$ is itself a presheaf. It is natural
to call it the universal presheaf though generally $P{\mcat T}$ is not
the classifying topos for $\mcat T$. It is easy to see that
$P{\mcat T}$ is a category with products (indeed limits)
defined pointwise. So we have the finite powers
${\mcat T}^m$ of the universal, with
${\mcat T}^m(n) = {\mcat T}(n)^m$.
The Yoneda Lemma is then the following.
\begin{proposition}\label{yoneda}
We have an isomorphism
$P{\mcat T}({\mcat T}^m, X) \cong X(m)$ natural in $X \in P{\mcat T}$.
\end{proposition}
For presheaves a Yoneda Lemma should lead to a Yoneda embedding.
This takes a particularly vivid form for the category $P{\mcat T}$.
The composition ${\mcat T}(n) \times {\mcat T}(m)^n \to {\mcat T}(m)$
corresponds by transpose to a map ${\mcat T}(n) \to 
({\mcat T}(m)^n \! \Rightarrow \! {\mcat T}(m))$
whose image is clearly natural in $m$. So we have
${\mcat T}(n) \to  P{\mcat T}({\mcat T}^n , {\mcat T})$ and
the Yoneda Lemma says in particular that this is an isomorphism.
But returning to the composition 
${\mcat T}(n) \times {\mcat T}(m)^n \to {\mcat T}(m)$,
we check directly that 
the diagram
\begin{diagram}
{\mcat T}(n) \times {\mcat T}(m)^n & \rTo &  {\mcat T}(m) \\
\dTo^{\cong} & & \dTo_{\cong} \\
P{\mcat T}({\mcat T}(n),{\mcat T}) \times 
P{\mcat T}({\mcat T}(m),{\mcat T})^n & \rTo & 
P{\mcat T}({\mcat T}(m),{\mcat T})
\end{diagram}
commutes. The result is the following form of an embedding
theorem.
\begin{proposition}\label{yonedaiso}
The Yoneda Lemma induces
an isomorphism between an algebraic theory $\mcat T$ and the endomorphism
theory of the universal object $\mcat T$ in $P{\mcat T}$.
\end{proposition}

Proposition \ref{yoneda} shows inter alia that the products of the 
universal are
dense, and so familiar arguments allow one to deduce more structure.
\begin{proposition} $P{\mcat T}$ is a topos; in
particular it is locally cartesian 
closed.
\end{proposition}
For the $\lambda$-calculus, we are interested in function spaces.
The following is essentially an old observation of Lawvere's,
the proof of which is an easy computation.
\begin{proposition}\label{lawvere}
For any presheaf $X$, the function space ${\mcat T}^p \Rightarrow X$
is given by the presheaf 
$({\mcat T}^p \Rightarrow X) (m) = X(m+p)$ with action 
$X(m+p) \times {\mcat T}(n)^m \to X(n+p)$
leaving the parameters $p$ undisturbed.  
\end{proposition}
The function spaces 
${\mcat T}^p \Rightarrow {\mcat T}$ with 
$({\mcat T}^p \Rightarrow {\mcat T})(m) = {\mcat T}(m+p)$
will be of particular interest for the $\lambda$-calculus.

\subsection{Monads and Lawvere Theories}\label{monlaw}
In this paper I use a multicategory theory approach
to algebraic theories. I regard it as the fundamental one
and it is particularly suited to the $\lambda$-calculus.
However for those who may prefer them I sketch the 
more traditional categorical approaches to algebra.

An algebraic theory $\mcat T$ induces a monad $T$ on $\bf Sets$ 
whose functor part is given by the coend formula 
$T(A) = \int^n {\mcat T}(n)
\times A^n$. 
Concretely $T(A)$ is the set of terms
from $\mcat T$ with constants from $A$ replacing variables.
The unit and multiplication of the monad correspond to constants
and substitution.
Conversely from a monad $T$ we get an algebraic theory
$\mcat T$ where ${\mcat T}(n) = T(n)$, since
terms in $n$ variables form the underlying set of the free
algebra on a set of size $n$. The monads which arise from 
algebraic theories are exactly the finitary monads
and there is an equivalence of categories between
algebraic theories and finitary monads. The
notions of algebra correspond. On the other
hand the presheaf topos which is fundamental to
my approach to the $\lambda$-calculus is less 
easy to handle from the monad point of view.

Lawvere theories introduced in \cite{Law63} provide
the other categorical approach to algebra. The equivalence 
with finitary monads was evident early but
Lawvere theories have been comparatively neglected.
That seems a mistake. \cite{HP07} compares the two approaches 
with applications to computer science in mind. With that
background, it is clear that the Lawvere Theory 
corresponding to an algebraic theory is essentially
(I suppress the identity on objects functor)
the free category $\cat T$ with products generated 
by $\mcat T$: one has simply ${\cat T}(n,m) = {\mcat T}(n)^m$
and the categorical composition is then easy.
Algebras are given by product-preserving 
functors ${\cat T} \to {\bf Sets}$ and almost by
definition this notion coincides with the one above.
For the presheaf categories there is the following easy consequence
of our Yoneda Lemma, Proposition \ref{yoneda}.
\begin{proposition}
The categories $P{\mcat T}$ 
and $P{\cat T} =
[{\cat T}^{\rm op},{\bf Sets}]$ are equivalent.
\end{proposition} 
Thus one could 
very readily rewrite this paper from the
Lawvere theory perspective.

\section{The lambda calculus}
\subsection{$\lambda$-theories}
The underlying philosophy is to take
seriously the basic syntax of term formation 
in context first made prominent in
Martin-L\"{o}f's treatment of Type Theory.
The term formation rules are of form
\[
 \dfrac{\Gamma \vdash t \in a \Rightarrow b
 \quad \Gamma \vdash u \in a}{\Gamma \vdash tu \in b}
 \hspace{1.5cm} 
\dfrac{\Gamma, \;  x \in a \vdash s \in b}
{\Gamma \vdash \lambda x.s \in a \Rightarrow b}
\]
while the computation rule is
$(\lambda x. s)u = s[u/x]$. The rules force one to 
handle terms in context and it is best to do that
directly. We can  capture the rules for $\lambda$-calculus
by saying that an interpretation
is a cartesian multicategory
with semi-closed structure. The pure lambda calculus 
corresponds to the one-object version.
\begin{definition}\label{lambda}
To equip an algebraic theory $\mcat L$ with {\em semi-closed structure}
is to give retractions ${\mcat L}(n+1) \, \ret \, {\mcat L}(n)$,
$\rho :{\mcat L}(n) \to {\mcat L}(n+1)$
and $\lambda: {\mcat L}(n+1) \to {\mcat L}(n)$ the {\em retraction}
and {\em section},
natural in $n$ and compatible with the actions
${\mcat L}(m) \times {\mcat L}(n)^m \to {\mcat L}(n)$ and
${\mcat L}(m+1) \times {\mcat L}(n)^m \to {\mcat L}(n+1)$.
A {\em $\lambda$-theory} is an
algebraic theory $\mcat L$ equipped with semi-closed
structure.
Let $\mcat L$ and $\mcat M$ be $\lambda$-theories. A {\em map
${\mcat L} \to {\mcat M}$ of $\lambda$-theories}
is a map of algebraic theories which commutes with retraction $\rho$
and section $\lambda$.
\end{definition}
We get a category of $\lambda$-theories. It is clearly
locally finitely presentable; but by the Fundamental Theorem \ref{fundthm}
it is in fact
equivalent to a category of algebras.

To start with one should understand the definition of $\lambda$-theory
concretely. The image of the identity ${\bf id} \in {\mcat L}(1)$ under 
$\rho: {\mcat L}(1) 
\to {\mcat L}(2)$ is a binary operation
${\bf app} \in {\mcat L}(2)$ of application, 
traditionally denoted by concatenation
${\bf app}(x,y) = xy$. 
Naturality implies that 
the retraction $\rho: {\mcat L}(n) \to {\mcat L}(n+1)$
is given by $a \to az$
where I adopt now and henceforth the convention to use $z$ as
the new extra variable in this and similar cases.
Starting with ${\bf id} \in {\mcat L}(1)$ and retracting $n$ times gives 
iterated applications ${\bf app}_n \in {\mcat L}(n+1)$ where we can write
suggestively ${\bf app}_{n+1} = {\bf app}({\bf app}_n,z)$ or spelling
out the variables 
${\bf app}_{n+1}(x, z_1, \cdots , z_{n+1}) = 
{\bf app}({\bf app}_n(x, z_1, \cdots , z_n), z_{n+1})$. By convention
the concatenation notation for application associates 
so that we could write $xz_1 \cdots z_n$.

The naturality of the section $\lambda: {\mcat L}(n + 1) \to {\mcat L}(n)$ 
for $\lambda$-abstraction is
more subtle. The critical fact is that 
${\mcat L}(2) \to {\mcat L}(1)$ does not generally take
application ${\bf app}$ to the identity ${\bf id}$. Rather 
it goes to an element $\app z$ where $\app \in {\mcat L}(0)$ 
is the image of ${\bf app}$ in ${\mcat L}(0)$. Naturality
then gives the following. Suppose $s \in {\mcat L}(n+1)$
has image $\lambda s \in {\mcat L}(n)$. Then we have $s = (\lambda s)z$
and $\lambda s = \app (\lambda s)$, and these conditions determine
$\lambda s$. So while maps of $\lambda$-theories
are not mere 
maps of the algebraic theories and the preservation
of the semi-closed structure is essential, 
we have the following.
\begin{proposition}\label{theorymap}
Suppose $\mcat L$ and $\mcat M$ are $\lambda$-theories.
If a map ${\mcat L} \to {\mcat M}$ of algebraic theories
preserves the binary operation ${\bf app}$ and constant $\app$,
then it is a map of $\lambda$-theories.
\end{proposition}
Of course one can iterate the section. The following will prove
essential. 
\begin{proposition}\label{abstraction}
Suppose that $s \in {\mcat L}(n)$ has image 
$\hat{s} = \lambda^n{s} \in {\mcat L}(0)$.
Then $s$ can be recovered from $\hat{s}$ and iterated
application via the identity
$s = {\bf app}_n (\lambda^n{s}, z_1, \cdots ,z_n)$.
\end{proposition}
Unsurprisingly the characterisation of $\lambda s$ extends. 
Rewriting the Proposition \ref{abstraction}, 
we have $s = (\lambda^n{s})z_1 \cdots z_n$,
and we also have $\app_n (\lambda^n{s}) = \lambda^n{s}$ where 
$\app_n = \lambda^{n+1}({\bf app}_n)$ is
the image in ${\mcat L}(0)$ of iterated application
$n$ arguments.
 
\subsection{Interpreting the $\lambda$-calculus}
I briefly sketch the interpretation of the lambda calculus 
in a $\lambda$-theory.
For clarity I shall here use semantic brackets $\lscott t \rscott$
to indicate the interpretation of a term $t$, though in accord with general
mathematical practice I shall dispense with them as soon as I prudently can.
Let $\mcat L$ be a $\lambda$-theory. Suppose that $\Gamma \vdash t$
is a term $t$ in context $\Gamma$ of length $n$. Then $t$ has an interpretation
$\lscott t \rscott \in {\mcat L}(n)$, which is defined inductively as follows.
For the $i$th variable 
$\Gamma \vdash x_i$, let $\lscott x_i \rscott = {\bf pr}_i$,
the $i$th projection. Application $tu$ of terms is interpreted by
the application so that 
$\lscott tu \rscott = {\bf app}(\lscott t \rscott, \lscott u \rscott )$;
and $\lambda$-abstraction of terms is defined using the section
$\lambda : {\mcat L}(n+1) \to {\mcat L}(n)$
so that $\lscott \lambda z.r \rscott
= \lambda(\lscott r \rscott)$.
Of course the point is not the definition but the 
fact that it works in the
sense that  $\beta$-equality is satisfied. 
The main point of the proof concerns substitution.
\begin{lemma} Suppose $\Gamma, z \vdash r$ and $\Gamma \vdash s$, so
we have interpretations $\lscott r \rscott \in {\mcat L}(n+1)$
and $\lscott s \rscott \in {\mcat L}(n)$. Then
$\Gamma \vdash r[s/z]$ has interpretation
$\lscott r[s/z] \rscott = 
\lscott r \rscott ({\bf pr}_1 , \cdots , {\bf pr}_n, \lscott s \rscott)$.
\end{lemma}
{\bf Proof}
By induction on the structure of $r$. The crucial $\lambda$-abstraction step
is just the compatibility of the sections 
$\lambda : {\mcat L}(n+1) \to {\mcat L}(n)$ with
the action.\\[0.4em]
There is a series of results which follow by easy
induction on the structure of $\lambda$-terms.
\begin{proposition}
The interpretation of the $\lambda$-calculus in a $\lambda$-theory
respects $\beta$-equality in the sense that 
$\lscott (\lambda z. s)u \rscott = \lscott s[u/z] \rscott$.
\end{proposition}
That provides moral support for the definition of $\lambda$-theory.
Next we consider maps of $\lambda$-theories.
\begin{proposition}
Suppose that $F: {\mcat L} \to {\mcat M}$ is a map of $\lambda$-theories.
Then $F$ preserves the interpretation of $\lambda$-terms:
for every $s$ we have 
$F\lscott s \rscott_{\mcat L} = \lscott s \rscott_{\mcat M}$.
\end{proposition}
Next we show that the 
$\lambda$-calculus presents the initial $\lambda$-theory. Let $\Lambda (n)$ 
be the terms of the $\lambda$-calculus in context of $n$ variables,
factored out by $\beta$-equality: write $[s] \in \Lambda(n)$ for 
the equivalence class
of $\Gamma \vdash s$, the terms $s$ in context $\Gamma$ of length $n$.
Identities and projections in $\Lambda$ are evident and composition is given
by substitution, so we have an algebraic 
theory. The retraction of $\Lambda(n)$ onto
$\Lambda (n+1)$ is $t \in \Lambda (n) \mapsto t.z \in \Lambda (n+1)$.
The section
is $s \in \Lambda (n+ 1) \mapsto \lambda z . s \in \Lambda (n)$.
We do indeed have a retraction because $(\lambda z.s)z = s$ in the
$\lambda$-calculus. Compatibility with the action is trivial
in the case of the retraction, while for the section it is the
basic syntactic lemma $(\lambda z.s)[r/x] = \lambda z (s[r/x])$.
So we get a $\lambda$-theory $\Lambda$. 
The last of our
results proved by structural induction is the following.
\begin{proposition} The interpretation $\lscott s \rscott \in \Lambda(n)$
of a $\lambda$-term $\Gamma \vdash s$ in context is given by its
equivalence class $[s]$.
\end{proposition}
From this essentially routine series of propositions
one gets the following basic result.
\begin{theorem}
The $\Lambda$ is the initial
$\lambda$-theory. 
\end{theorem}
Henceforth I shall use standard $\lambda$-calculus notation to define
elements
in $\lambda$-theories. We have already seen $\app = \lambda xy. xy$ (and
more generally $\app_n = \lambda z_1 \cdots z_{n+1}.z_1 \cdots z_{n+1}$).

Note that the definition of $\lambda$-theory involves an unfamiliar 
use of the notion of
algebraic theory. The initial algebraic 
theory is the pure 
theory of equality
and usually we extend that with constants and function symbols and equation
between them. Here nothing of that kind is involved yet what we get
is highly non-trivial. It is perhaps just worth drawing attention to a tension
with usual examples of algebra. Operations
like $x \mapsto x+1$ generally have no fixed points. But the $\lambda$-calculus
has fixed point operators. So for example the only semi-closed theory
of rings is the terminal theory with $0=1$.

\subsection{Extensionality}

I turn aside for a moment to
comment on two separate issues termed extensionality.
What is called weak extensionality caused much concern
in early developments; note the faintly apologetic remark
`In spite of not being weakly extensional $\lambda$-algebras
are worth studying' on page 87 of \cite{Bar84}.
The category theoretic reading of weak extensionality
is the question whether one has enough points
\cite{Scott80}, and it is the same for $\lambda$-theories.
An algebraic theory $\mcat T$ {\em has enough points} 
just when equality on each ${\mcat T}(n)$ is reflected in
the action
${\mcat T}(n) \times {\mcat T}(0)^n \to {\mcat T}(0)$
on constants. $\lambda$-theories with enough points are
essentially the Frege Structures of \cite{Aczel80}
but $\lambda$-theories may not have enough
points: indeed the initial $\lambda$-theory does not. 
One does not need to worry about that. For those who do,
there is this general fact.
\begin{proposition}
Any algebraic theory $\mcat T$ or $\lambda$-theory
$\mcat L$ embeds in an algebraic theory or $\lambda$-theory 
with enough points. 
\end{proposition}
{\bf Proof} Take 
${\mcat T}(\omega)$ to be the free algebra on countable 
many generators, constructed in the obvious
way as a direct limit of free algebras ${\mcat T}(n)$.
Then $\mcat T$ embeds in ${\mcat T}_{{\mcat T}(\omega)}$,
the theory of extensions of ${\mcat T}(\omega)$, which
clearly has enough points. A $\lambda$-theory $\mcat L$ embeds 
in ${\mcat L}_{{\mcat L}(\omega)}$, and the $\lambda$-theory
structure is easy.\\[0.4em]
Of course this observation is not new for the $\lambda$-calculus: it
is essentially in Section 4 of \cite{BK80}.
A special case is the move from the initial closed term to the open term 
interpretation.

Another quite different aspect of extensionality
is the $\eta$-rule $\lambda x. tx = t$
when $x$ is not free in $t$. This
corresponds exactly to the requirement that the semi-closed 
structure be closed, in the sense that the retractions
$\rho: {\mcat L}(n) \to {\mcat L}(n+1)$ are isomorphisms.
In terms of the analysis above this is equivalent to 
application ${\bf app}$ being taken to the identity $\bf id$.
The theory outlined in this paper is not much affected by adding
this condition but doing so confuses
rather than helps. There can be no doubt that it is the basic
calculus with just the $\beta$-rule which is fundamental. 

\subsection{Reflexive objects} 
Scott's elegant categorical approach to the semantics of the
$\lambda$-calculus stands out as clean mathematics. It
is not flexible as a definition,
but it captures much of importance.
The most obvious element is a Representation
Theorem. I give a new proof: the $\lambda$-theory
perspective makes Scott's idea easy to understand. 
Suppose that $\cat C$ is a 
cartesian closed category and $U$
an object of $\cat C$ equipped with a retraction onto
the function space $U^U = U \! \Rightarrow \! U$. Set 
${\mcat U}(n) = {\cat C}(U^n,U)$. This is automatically
an algebraic theory, the endomorphism
theory of $U$. Moreover since 
${\cat C}(U^n, U^U)  \cong {\cat C}(U^{n+1}, U)$ we get
retractions
${\mcat U}(n+1) \ret {\mcat U}(n)$ manifestly natural in $n$.
So we get a $\lambda$-theory $\mcat U$, the {\em endomorphism
$\lambda$-theory} of the reflexive object $U$.

The essence of
Scott's wonderful insight was that any $\lambda$-theory 
can be so represented. His proof with discussion
of the significance of the result is in \cite{Scott80}. 
It involves serious coding which we shall not need until Section \ref{monoid}.
\begin{theorem} {\rm (Scott's Representation Theorem)} 
Any $\lambda$-theory is isomorphic to the endomorphism
$\lambda$-theory of a reflexive object in
a cartesian closed category.\label{rep1}
\end{theorem}
{\bf Proof} We prepared for this in
Section \ref{topos}. Given $\mcat L$ take 
the presheaf topos $P{\mcat L}$ and let $U$ be the
universal object $\mcat L$ itself. Proposition \ref{lawvere}
gives the function space
$U^U$: it is the presheaf ${\mcat L}(n+1)$. By
definition a $\lambda$-theory consists of a retraction of
${\mcat L}(n)$ onto ${\mcat L}(n+1)$ and the naturality
conditions say that the maps are in
$P{\mcat L}$. So we have a retraction of $U$ onto $U^U$
and $U$ is a reflexive object.
It remains to consider the $\lambda$-theory $\mcat U$
obtained from $U$. We have 
${\mcat U}(n) = {P{\mcat L}}(U^n,U) \cong {\mcat L}(n)$,
an isomorphism of algebraic theories 
by Proposition \ref{yonedaiso}, and then evidently
an isomorphism of $\lambda$-algebras.\\[0.4em]
The proof \cite{Scott80}
in the Curry Festschrift gives the reflexive object
in the cartesian closed category of retracts of a $\lambda$-algebra.
For a $\lambda$-theory version, take the monoid ${\mcat L}(1)$ 
as a one object category. Any category $\cat C$ has
a category of retracts, its Karoubi envelope or Cauchy completion.
The objects are idempotents $e: A \to A$ with $e \circ e = e$ in $\cat C$;
and maps between idempotents $e: A \to A$ 
and $f: B \to B$ are maps $v : A \to B$ such that $f \circ v \circ e = v$.
Composition is inherited from $\cat C$ and each idempotent $e$
is its own identity. Scott's category $\bf R$ of retracts
is the Cauchy completion of
${{\mcat L}(1)}$. Scott showed by explicit calculation
that it is cartesian closed. I shall use the presheaf category
$P{\mcat L}$ to give a new proof of this, see Corollary \ref{retracts}.

\subsection{The Taylor Fibration}\label{taylor}

In his PhD thesis \cite{Tay86},
Paul Taylor extended Scott's analysis in a very remarkable way. 
Taylor observed that the category ${\cat R}$ of retracts 
is not just cartesian closed but is relatively cartesian closed. 
I briefly explain what that means. Taylor gave a notion of a 
category $\cal C$ {\em equipped with display maps}
$\cal D$ with the idea that the display maps model dependent or indexed types
for some type theory. This approach is now standard: under the influence of
homotopy type theory display maps are often now called fibrations but I shall
stick with Taylor's terminology. Basic properties of a type theory
are given by closure properties. Taylor required that $\cal D$ be closed by 
pullback along all maps, so indexed types are closed under 
arbitrary substitution; that $\cal D$ be closed under composition,
so that types are closed under indexed sums; and finally that $\cal D$ 
contain all terminal projections, a condition which is best thought
of as avoiding redundancy. With that in place, Taylor defined a
category $\cal C$ with displays $\cal D$ to be {\em relatively
cartesian closed} if the display maps are closed by products along
display maps.

For each object $E \in {\cat R}$,
Taylor localised the construction of $\cat R$.
Let $U \in {\cat R}$ be the generating object given by the
identity idempotent. Then over each $E \in {\cat R}$ we have
$\Delta_E  (U) = (U \times E \to E)$ which is a reflexive 
object in the slice ${\cat R}/E$. Taylor considered 
the subfibration ${\cat R}(E)$
of ${\cat R}/E$ consisting of retracts $A \to E$ of $\Delta_E  (U)$.
His result is that for every $\alpha: F \to E$ in ${\cat R} (E)$
the pullback functor $\alpha^* :{\cat R}(E) \to {\cat R}(F)$, which is
evident from the definition of ${\cat R}(E)$, comes equipped
with a right adjoint  $\Pi_{\alpha} : {\cat R}(F) \to {\cat R}(E)$
and left adjoint 
$\Sigma_{\alpha} : {\cat R}(F) \to {\cat R}(E)$. That
$\cat R$ is cartesian closed is a simple consequence.

In the spirit of \cite{Scott80}, Taylor
simply wrote down the various
combinators and calculated to show that they work. 
I shall give a more abstract proof using the presheaf category
$P({\mcat L})$ with universal object $U$. 
We already saw the retract from $U$ to $U^U$. We need also
a retract from $U$ to $U \times U$.
It relies on familiar $\lambda$-calculus: we have the section
\[
{\mcat L}(n) \times {\mcat L}(n) \to {\mcat L}(n) ; \quad
(a,b) \to \lambda x. xab \, ,
\]
and the retraction
\[
{\mcat L}(n) \to {\mcat L}(n) \times {\mcat L}(n) ; \qquad
c \to (c {\bf T}, c {\bf F}) \, ,
\]
where ${\bf T} = \lambda xy . x$ and ${\bf F} = \lambda xy . y$.

Now consider Taylor's fibration over the whole of $P({\mcat L})$.
For $X \in P({\mcat L})$, take ${\cat R}(X)$ to be the 
category of retracts
of the reflexive object $\Delta_X(U) = (U \times X \to X)$ in 
$P({\mcat L})/X$.
\begin{theorem} Take $X \in P({\mcat L})$. 
Let $\alpha : Y \to X$ be in ${\cat R}(X)$. Then for any $Q \to Y$ in
${\cat R}(Y)$, the standard indexed product $\Pi_{\alpha} Q \to X$ and 
standard indexed sum $\Sigma_{\alpha} Q \to X$ over 
$X \in P({\mcat L})$ both lie in ${\cat R}(X)$.
\end{theorem}
{\bf Proof} The aim is to reduce to a special case.
First $Q \to Y$ is a retract of $\Delta_Y (U)$, so
$\Pi_{\alpha} Q \to X$ and $\Sigma_{\alpha} Q \to X$ are
retracts of $\Pi_{\alpha} \Delta_Y (U)$ and 
$\Sigma_{\alpha} \Delta_Y (U)$, respectively; so it suffices 
to show that these lie in ${\cat R}(X)$. Now consider a diagram
\begin{diagram}
Y & \rTo^{i} & U \times X & \rTo^{r} & Y \\
 & \rdTo_{\alpha} & \dTo^{q} &\ldTo_{\alpha} & \\
 && X &&
\end{diagram}
displaying $\alpha: Y \to X$ as a retract of $\Delta_XU$. From it we get maps
\[
\Pi_{\alpha} \Delta_Y U \to \Pi_{\alpha}\Pi_r \Delta_r \Delta_Y U
\cong \Pi_q \Delta_{U \times X}U \to 
\Pi_q \Pi_i \Delta_i \Delta_{U \times X}U \cong \Pi_{\alpha} \Delta_Y U 
\]
\[
\Sigma_{\alpha} \Delta_Y U \cong 
\Sigma_q \Sigma_i \Delta_i \Delta_{U \times X} U \to 
\Sigma_q \Delta_{U \times X} U \cong 
\Sigma_{\alpha} \Sigma_r \Delta_r \Delta_Y U \to
\Sigma_{\alpha} \Delta_Y U
\]
displaying retracts 
$\Pi_{\alpha} \Delta_Y U \, \ret \, \Pi_q \Delta_{U \times X}U$ and 
$\Sigma_{\alpha} \Delta_Y U \, \ret \, \Sigma_q \Delta_{U \times X} U$.
So it suffices for products to consider $\Pi_q \Delta_{U \times X}U 
\cong \Delta_X(U \Rightarrow U)$; but we know $(U \Rightarrow U) \ret U$,
and so $\Delta_X(U \Rightarrow U) \ret \Delta_X U$. 
Similarly for sums we consider $\Sigma_q \Delta_{U \times X}U 
\cong \Delta_X(U \times U)$; but we saw above that $U \times U \ret U$,
and so $\Delta_X(U \times U) \ret \Delta_X U$. That completes the proof.

\begin{corollary} {\rm (Paul Taylor)} The category of retracts of a 
$\lambda$-theory is relatively cartesian closed.
\end{corollary}
{\bf Proof} The category $\cat R$ of retracts of $U$
is a subcategory of $P({\mcat L})$. Furthermore if 
$X \in {\cat R}$
then so are the objects of 
${\cat R}(X)$. So the result is immediate by restricting to $\cat R$.

\begin{corollary}\label{retracts} {\rm (Dana Scott})
The category of retracts of a $\lambda$-theory is cartesian closed. 
\end{corollary}
{\bf Proof} Immediate by restricting 
to the fibre ${\cat R} = {\cat R}(1)$ over $1$.

In Section \ref{monlaw}, I remarked that for a general algebraic
theory $\mcat T$, the presheaf category $P{\mcat T}$ is equivalent
to the standard presheaf category $P{\bf T}$ on the Lawvere theory
$\bf T$ generated by $\mcat T$.
For $\lambda$-theories $\mcat L$ we have much more.
Since the category $\bf R$ of retracts is closed under products
we have a functor ${\bf L} \to {\bf R}$ from the Lawvere theory
$\bf L$ generated by $\mcat L$. Moreover the monoid
${\mcat L}(1)$ is embedded in the $\lambda$-theory $\mcat L$.
Restriction gives functors 
$P{\bf R} \to P{\bf L} \to P{\mcat L} \to P{\mcat L}(1)$. 
\begin{proposition}\label{morita}
The functors $P{\bf R} \to P{\bf L} \to P{\mcat L} \to P{\mcat L}(1)$
are equivalencs.
\end{proposition} 
{\bf Proof} ${\mcat L}(1)$ embeds in $\bf L$
which embeds in $\bf R$ the category of retracts of ${\mcat L}(1)$,
so $P{\bf R} \simeq P{\bf L} \simeq P{\mcat L}(1)$ by Morita theory. 
The equivalence $P{\mcat L} \simeq P{\cat L}$ from Section \ref{monlaw} is
trivial.\\[0.4em]
I shall make use of the equivalence $P{\mcat L} \to P{\mcat L}(1)$
at the beginning of Section \ref{fundsec}. It takes the explicit
form
\[
\big( X(n) \times {\mcat L}(m)^n \to X(m) \big) \; \mapsto \;
\big( X(1) \times {\mcat L}(1) \to X(1) \big) \, . 
\]

\section{Algebras}

\subsection{Algebras for $\lambda$-theories}

Consider any $\lambda$-theory $\mcat L$ simply 
as an algebraic theory,
and then we have a category ${\rm Alg}({\mcat L})$
of $\mcat L$-algebras. 
A $\mcat L$-algebra is a set $A$ equipped with actions 
${\mcat L} (n) \times A^n \to A$. 
Concretely that means that for each term $t({\bf x})$ 
in ${\mcat L} (n)$ and each $n$-tuple ${\bf a} \in A^n$
we get an interpretation $t({\bf a})$ in $A$; and this
behaves as expected on variables and respects substitution
and $\beta$-equality.
It turns out that one can focus almost entirely
on  $\Lambda$, the initial $\lambda$-theory.
\begin{definition}
A {\em $\Lambda$-algebra} is an algebra for the the initial
$\lambda$-theory $\Lambda$.
\end{definition}
This is a clean definition close in spirit to environment
or valuation models but avoiding the explicit interpretation
of abstraction and so the issue of weak 
extensionality. Giving $\Lambda$-algebra structure 
amounts to giving $\lambda$-algebra structure as in \cite{Bar81},
but to underline the difference in perspective, I use the category theoretic
terminology. Map of $\Lambda$-algebras are 
maps in the category ${\rm Alg}(\Lambda)$. This is
so smooth that one might wonder why one should not take this
as the basic definition. The answer is that one never
directly shows that one has a $\Lambda$-algebra. One needs an
induction over $\lambda$-terms with free variables and 
that amounts to considering $\lambda$-theories.

I say a bit more about what the definition means concretely.
Think of $s \in \Lambda (n)$ as a $\lambda$-term with $n$
free variables. A $\Lambda$-algebra $A$ is equipped with an
interpretation $s({\bf a}) \in A$ of the $\lambda$-term 
for every term $s$ with constant ${\bf a} \in A^n$ substituted
for the free variables. In particular there is an interpretation
$(a,b) \mapsto ab$ of the application $\bf app$ as a binary operation
and also interpretations of all constant $s \in \Lambda (0)$
as $s \in A$. That determines the structure in the sense of the
following.
\begin{lemma}\label{appconst}
Suppose that $A$ and $B$ are $\Lambda$-algebras and $f:A \to B$
a map preserving the binary operation of application and the 
$\lambda$-definable constants. Then $f$ is a map of
$\Lambda$-algebras. 
\end{lemma} 
{\bf Proof} This comes from Proposition \ref{abstraction}.
Write $s({\bf x})$ as ${\bf app}_{n+1}(\hat{s}, {\bf x})$
where $\hat{s} = \lambda^n s$ is the constant obtained by
iterated 
$\lambda$-abstraction. The interpretation of constants and 
iterated application are preserved by $f$. Hence so is the
interpretation of $s$.\\[0.2em]

Any $\lambda$-theory $\mcat L$ gives rise to a $\Lambda$-algebra
in a straightforward way. ${\mcat L}(0)$ is the initial $\mcat L$-algebra
and composition with the unique $\Lambda \to {\mcat L}$ makes it a
$\Lambda$-algebra. Here it is only the map $\Lambda \to {\mcat L}$
which matters and ${\mcat L}$ need not be a $\lambda$-theory. But that
generality is of no significance. The following is trivial.
\begin{proposition} 
The operation ${\mcat L} \mapsto {\mcat L}(0)$ gives a functor
from $\lambda$-theories to $\Lambda$-algebras.
\end{proposition}

\subsection{Presheaves on the monoid}\label{monoid}

Let $A$ be a $\Lambda$-algebra. On 
${A}(1) = \{ a\in A \, | \, \app a = a \}$ take the monoid structure
with multiplication
$(a,b) \mapsto a \circ b = \lambda x. a(bx)$
representing composition. Write $M_A$ for this monoid.
The underlying set is a 
retract of $A$ and another way to give $M_A$
is as formal elements of the form $az$ factored out by
$a \sim b$ if and only if $\app a = \app b$ and with
composition $(az, bz) \mapsto a(bz) = (a \circ b)z$.
Thus for a $\lambda$-theory $\cal L$, the monoid 
$M_{{\mcat L}(0)}$ is isomorphic to the monoid ${\mcat L}(1)$. 
Motivated by Section \ref{taylor},
consider $PA = PM_A$ the category of presheaves
on the one object category $M_A$.
The universal object $U = U_A$ is $A(1)=M_A$ as underlying set 
with the evident right action of $b \in M_A$ by composition, 
$(a,b) \mapsto a \circ b$. 

The association $A \mapsto PA$ extends to maps. 
Suppose $f: A \to B$ is a map of 
$\Lambda$-algebras. As $f$ preserves application and $\app$,
it maps $A(1)$ to $B(1)$. Moreover since composition is implemented by
$\lambda xyz. x(yz)$, we get a map
$M_f : M_A \to M_B$ of monoids. Now composition with $M_f$
gives an easy functor $PB = PM_B \to PM_A = PA$ and 
this has a left adjoint $Pf:PA \to PB$ given by left Kan extension.
I give the following now by way of a warm-up:
some of the calculations appear in \cite{Koy84}.
\begin{proposition}\label{prodpres}
For a map of $\Lambda$-algebras $f: A \to B$, the induced $Pf: PA \to PB$
preserves finite products.
\end{proposition}
{\bf Proof}
It is sufficient to show preservation of
the terminal object and
preservation of products of representables. 
These
correspond to conditions on categories of elements
as follows.\\
{\em Terminal object}. Consider the category with objects
$b \in M_B$ and maps $b \to \bar{b}$ being given by $a \in M_A$
such that $f(a)\circ b = \bar{b}$. We need to show that this is
connected. But for every $b \in M_B$, $\lambda x. \cI \in M_A$ 
gives a map from $b$ to $\lambda x. \cI \in M_B$, so the 
latter is weakly terminal. (I use $\cI = \lambda x.x$, the
identity combinator but any $\lambda$-definable 
constant will do.)\\
{\em Binary Products}. For $b_1,b_2 \in M_B$ consider the
category with objects given by $c \in M_B$ and
$a_1,a_2 \in M_A$ with $f(a_1) \circ c = b_1$ and
$f(a_2) \circ c = b_2$. Maps  from $c \in M_B$ and
$a_1,a_2 \in M_A$ to  $\bar{c} \in M_B$ and
$\bar{a}_1,\bar{a}_2 \in M_A$ are given by elements $\hat{a}$ in
$M_A$ with $f(\hat{a}) \circ c = \bar{c}$ and. We need to show this 
category is connected.
Take $s= \lambda w(\lambda x. x (b_1 w) (b_2 w))$ in $M_B$
and $p = \lambda x.x{\bf T}$ and $q = \lambda x . x{\bf F}$
in $M_A$. Since
\[
p \circ s =
\lambda w.p(\lambda x. x (b_1 w) (b_2 w)) =
\lambda w((\lambda x. x (b_1 w) (b_2 w)){\bf T}) =
\lambda w. b_1 w = b_1 \, ,
\]
and similarly $q \circ s = b_2$, 
we have an object of our category. Moreover given any object $c \in M_B$ and
$a_1,a_2 \in M_A$ as above we can take 
$r = \lambda w(\lambda x. x (a_1 w) (a_2 w))$. Similar calculations 
to those just given show $p \circ r = a_1$, $q \circ r = a_2$ and
finally one checks $f(r) \circ c = s$. So again we have found
a weakly terminal object.

\subsection{The function space analysis}\label{fnspace}
For $A$ a $\Lambda$-algebra, we study the function space 
$U^U$ of the universal $U \in PA$. 

First I give a general categorical analysis. Let $M$ be a monoid
and $X,Y \in PM$ presheaves on $M$. The function space
$Y^X$ can be represented as the set $PM(M \times X, Y)$ of
$M$-equivariant maps, that is of $\phi:M \times X \to Y$
such that $\phi(m.m',x.m') = \phi(m,x).m'$, with action given
by $\phi . \bar{m} (m,x) = \phi(\bar{m}.m,x)$. This makes sense
since
\[
(\phi.\bar{m}(m,x)).m' = \phi(\bar{m}.m.m',x.m') = \phi.\bar{m}(m.m',x.m')
\]
and so $\phi.\bar{m}$ is again $M$-equivariant. Finally the evaluation
map $Y^X \times X \to Y$ is given in this representation by
\[
PM(M \times X , Y ) \times X \to Y \, ; \quad (\phi , x) \mapsto 
\phi(I,x)
\]
where here $I$ is the unit of the monoid $M$.

For the monoid $M_A$ of
a $\lambda$-algebra $A$ there is a more concrete
representation of the function space $U^U$.
Let $A(2) = \{ d \in A \, | \, {\app}_2 d=d \} = 
\{d \in M_A \, | \, \app \circ d =d \}$. The second
characterisation shows that $A(2)$ has an action of $M_A$ by 
composition on the right.
Now given $d \in A(2)$ we define a corresponding 
$\phi : M_A \times M_A \to M_A$
by $\phi(a,b) = \lambda y. d(ay)(by)$. Clearly
\[
\phi(a\circ c, b \circ c) = \lambda y. d(a(cy))(b(cy)) = \phi(a,b) \circ c
\]
so $\phi$ is equivariant. Thus $d \mapsto \phi$ is a map of sets 
$A(2) \to PA(U \times U, U)$. Furthermore for $c \in M_A$, 
$d \circ c = \lambda xy.d(cx)y$ maps to the function
$\lambda y.d(c(ay))(by) = \phi.c(a,b)$ of $a,b$. Thus
$d \mapsto \phi$ gives a map in $PA$ from $A(2)$ to $U^U$.

Now let $p = \lambda x. x{\bf T}$ and $q = \lambda x. x{\bf F}$. Reflecting on
the argument of Section \ref{taylor} leads one to think that 
$\phi(p,q)$ is in some sense generic.
I exploit that thought.
\begin{proposition} The map $A(2) \to U^U$ above is an isomorphism
in $PA$.
\end{proposition}
{\bf Proof}
Given a map $\phi : A(1) \times A(1) \to A(1)$, set
$d = \lambda yz. \phi(p,q)(\lambda x. xyz) \in A(2)$. This provides an
inverse to our map $A(2) \to U^U$ in $PA$. 
For first starting with $d$ passing to $\phi$
and back gives
\begin{equation*}
\begin{split}
\lambda yz.\big( \lambda w. d (pw)(qw)\big) (\lambda x.xyz) & =
\lambda yz.\big( \lambda w. d (w{\bf T})(w{\bf F})\big)(\lambda x.xyz) \\
& = \lambda yz. d({\bf T}yz)({\bf F}yz) \\
& = \lambda yz. dyz = d \\ 
\end{split}
\end{equation*}
On the other hand starting with $\phi$ passing to $d$ and back
gives a map taking $(a,b)$ to
\begin{equation*}
\begin{split}
\lambda w.\big( \lambda yz.(\phi (p,q)(\lambda x.xyz)\big) (aw)(bw)& =
\lambda w.\phi (p,q)(\lambda x.x(aw)(bw))\\
& = \phi (p,q) \circ \lambda w.(\lambda x.x(aw)(bw))\\
& =  \phi \big( p \circ \lambda w.(\lambda x.x(aw)(bw)),q \circ 
\lambda w.(\lambda x.x(aw)(bw)) \big) \\
& =  \phi \big( \lambda w.p(\lambda x.x(aw)(bw)), 
\lambda w.q(\lambda x.x(aw)(bw)) \big) \\
& = \phi(a,b)\\
\end{split}
\end{equation*}
as
$\lambda w.p(\lambda x.x(aw)(bw))
= \lambda w.{\bf T}(aw)(bw) = \lambda w.aw = a$,
and
$\lambda w.q(\lambda x.x(aw)(bw)) = b$ similarly.
So our map is bijective on underlying sets and so an isomorphism.\\[0.3em]
Again we see calculations from \cite{Koy84}.
We shall also need the precise form of evaluation arising from the 
identification of $A(2)$ with $U^U$. It is
\[
A(2) \times A(1) \to A(1) \, ; \quad (d,a) \mapsto 
\lambda y. d({\bf I}y)(ay) =
\lambda y. dy(ay) \, ,
\]
as the identity combinator ${\bf I} = \lambda x.x$ is the unit
in $M_A$.

For the record there is an easy extension of
the analysis of $U^U$. Each $U^n \Rightarrow U$ can be represented by
${A}(n+1) = \{ d | \app_n \circ d = d \}$ with again the obvious
action. Under evaluation $d \in A(n+1)$ corresponds to
$(a_1, \cdots , a_n) \mapsto \lambda y. dy(a_1y) \cdots (a_ny)$
as a map $U^n \to U$. Observe that generally 
$\{ d \in M_A | \app_n \circ d = d \} = \{d \in A | \app_{n+1}d = d \}$,
where the left hand side has a clear action on the right by
the monad $M_A$.

\subsection{The $\lambda$-theory of a $\Lambda$-algebra}
In the previous section we saw that the function space $U^U$ of
the universal object $U$ in $PA$ is given 
by $A(2) = \{d \in M_A | \app \circ d = d \}$ with the 
action of $M_A$ on the right. Evidently composition 
on the left with $\app$ gives a retract from $U$ to $U^U$
and the generic $U$ is a reflexive object in the 
presheaf category $PA$. 

\begin{definition}
The $\lambda$-theory ${\mcat U}_A$ 
of a $\Lambda$-algebra $A$ is the theory 
of the reflexive universal $U \in P(A)$. 
\end{definition}
Take $f: A \to B$ a map of $\Lambda$-algebras
with induced functor $Pf: PA \to PB$.
The left Kan extension $Pf$ takes the universal $U_A$ in $P(A)$ to 
$U_B$ in $P(B)$ with a specified isomorphism
$Pf(U_A) \cong U_B$. So
$A \mapsto (PA,U_A)$ is pseudofunctorial with $(PA,U_A)$
considered as a category with specified object. 
By Proposition \ref{prodpres} $Pf$ preserves finite
products and so gives maps 
$P(A)(U_A^n,U_A) \to P(B)(U_B^n,U_B)$
which taken together give a map of algebraic 
theories ${\mcat U}_A \to {\mcat U}_B$.
\begin{proposition}\label{universal}
The operation $A \mapsto {\mcat U}_A$ gives a functor
from $\Lambda$-algebras to $\lambda$-theories. 
\end{proposition}
{\bf Proof}
Given a map $f:A \to B$ of $\Lambda$-algebras, we first 
check that the induced map ${\mcat U}_A \to {\mcat U}_B$
is a map of $\lambda$-theories. But $f$ preserves $\app$
which determines the function space as a retract of the
universal. So $Pf$ preserves the retract $U^U \ret U$
and the result follows. Furthermore $A \mapsto M_A$ is functorial
in $A$, so $A \mapsto (PA,U_A)$ is pseudofunctorial, but then
$A \mapsto {\mcat U}_A$ is functorial as we have a mere
category of algebraic theories.

\subsection{The Fundamental Theorem}\label{fundsec}

Let ${\mcat L}$ be a $\lambda$-theory. Composing the Yoneda with
the equivalence $P{\mcat L} \to P{\mcat L}(1)$ of Section \ref{taylor}
gives isomorphisms
${\mcat L}(n) \to P{\mcat L}(1)(U^n,U)$ taking
$a \in {\mcat L}(n)$ to the map
${\mcat L}(1)^n \to {\mcat L}(1)$ given by 
$(b_1, \cdots , b_n) \mapsto a(b_1, \cdots , b_n)$.
This gives an isomorphism between $\mcat L$ and the
endomorphism $\lambda$-theory of $U \in P{\mcat L}(1)$.
Furthermore we have a canonical isomorphism
$M_{{\mcat L}(0)} \cong {\mcat L}(1)$ of monads and so
an isomorphism $P{\mcat L}(1) \cong P{\mcat L}(0)$.
Thus we get isomorphisms
${\mcat L}(n) \to P{\mcat L}(0)(U^n,U) = {\mcat U}_{{\mcat L}(0)}(n)$
taking
$a \in {\mcat L}(n)$ to the map
${\mcat L}(0)^n \to {\mcat L}(1)$ given by 
$(c_1, \cdots , c_n) \mapsto \lambda x. a(c_1x, \cdots , c_nx)$.
This gives an isomorphism of $\lambda$-theories
$\eta_{\mcat L} :{\mcat L} \to {\mcat U}_{{\mcat L}(0)}$. The following
is immediate.
\begin{proposition}\label{eta}
The $\eta_{\mcat L} :{\mcat L} \to {\mcat U}_{{\mcat L}(0)}$
are $\lambda$-theory isomorphisms natural in 
$\mcat L$.
\end{proposition}

Now take a 
$\Lambda$-algebra $A$, pass to the universal $\lambda$-theory
${\mcat U}_A$ and then take the induced $\Lambda$-algebra
${\mcat U}_A(0)$. \cite{Scott80} notes that
these are isomorphic, and details of a syntactic argument
are given in \cite{Bar84} and \cite{Koy84}.
One considers open terms in the $\lambda$-calculus
with constants from $A$ and shows inductively that each
term $t({\bf x})$ with $n$ free variables is interpreted
in ${\mcat U}_A(n) \cong A(n)$ as above by the interpretation
in $A$ of its closure $\lambda {\bf x}. t({\bf x})$.
The spirit of categorical logic is to engage in as few syntactic
inductions as possible, and I present an alternative.
The $\Lambda$-algebra 
${\mcat U}_A(0) = PA(1,U)$ consists of the fixed points 
of $A(1)$ under the composition action. There is evidently a map
of sets $\varepsilon_A : {\mcat U}_A(0) \to A; a \mapsto a{\bf I}$.
(Since $a$ is fixed any constant will do.)
\begin{lemma}\label{application}
The maps $\varepsilon_A : {\mcat U}_A(0) \to A $ preserve the application.
\end{lemma}
{\bf Proof}
Internal application $U \times U \to U$ is
given by evaluation. It follows that it takes $(a,b) \in A(1)^2$
to $\lambda y. ay(by)$. But $\varepsilon_A (\lambda y. ay(by)) =
a{\bf I}(b{\bf I}) = \varepsilon_A(a)\varepsilon_A(b)$.\\[0.3em]
We do not yet know that the $\varepsilon_A$ are maps of $\Lambda$-algebras,
but for $f :A \to B$
a map of $\Lambda$-algebras the naturality
diagram
\begin{diagram}
{\mcat U}_A(0)& \rTo^{{\mcat U}_f(0)} & {\mcat U}_B(0) \\
\dTo^{\varepsilon_A} & & \dTo^{\varepsilon_B}\\
A & \rTo_f & B
\end{diagram}
commutes because $f(a{\bf I}) = f(a){\bf I}$. It is an equally
trivial calculation to show that for $\mcat L$ a $\lambda$-theory
the familiar triangle identity diagram 
\begin{diagram}
{\mcat L}(0) & \rTo^{\eta_{L}(0)} & {\mcat U}_{{\mcat L}(0)}(0) \\
& \rdEq & \dTo_{\varepsilon_{{\mcat L}(0)}}\\
& & {\mcat L}(0)
\end{diagram}
commutes: specifically the composite is 
$a \mapsto \lambda x.a \mapsto (\lambda x.a){\bf I} = a$.
\begin{proposition}\label{eps}
The $\varepsilon_A :{\mcat U}_A(0) \to A$ are $\Lambda$-algebra 
isomorphisms natural in $A$.
\end{proposition}
{\bf Proof} To show that
$\varepsilon_A$ is a map of $\Lambda$-algebras, consider the commuting diagram
\begin{diagram}
\Lambda(0) & \rTo^{\eta_{\Lambda}(0)} 
& {\mcat U}_{\Lambda(0)}(0) & \rTo &{\mcat U}_A(0) \\
 & \rdEq & \dTo_{\varepsilon_{\Lambda(0)}} & & \dTo_{\varepsilon_A} \\
& & \Lambda(0) & \rTo & A \, .
\end{diagram}
We have unique $\Lambda$-algebra maps $\Lambda (0) \to {\mcat U}_A(0)$
and $\Lambda (0) \to A$. It follows that $\varepsilon_A$ preserves
$\lambda$-definable constants. By Lemma \ref{application} it
preserves application, so by Proposition
\ref{appconst} it is indeed a $\Lambda$-algebra map. We saw the trivial
naturality above and $a \mapsto a{\bf I}$ has inverse $c \mapsto \lambda x. c$
so we have a natural isomorphism.\\[0.2em]

It follows from Propositions \ref{eta} and \ref{eps} that the functors 
${\mcat L} \mapsto {\mcat L}(0)$ and
$A \mapsto {\mcat U}_A$ give an equivalence (in fact an adjoint
equivalence)
between the categories of $\lambda$-theories and of $\Lambda$-algebras.
That is almost but not quite what I want to 
call the Fundamental Theorem: there is a little bit more.
Consider $\Lambda_A$ the algebraic theory of extensions of $A$
as described in Section \ref{algebras}.
It is not a priori obvious that it is a $\lambda$-theory
but we can identify it with ${\mcat U}_A$. There are a number of ways to
see this in terms of the equivalence. The following seems
down to earth.
For any $\Lambda$-algebra $A$ we have
the unique map of $\lambda$-theories $\Lambda \to {\mcat U}_A$.
Using the isomorphism $A \cong {\mcat U}_A(0)$ we get a factorization
$\Lambda \to {\Lambda}_A \to {\mcat U}_A$ of algebraic theories. We get an
induced functor ${\rm Alg}({\mcat U}_A) \to{\rm Alg}({\Lambda}_A)$.
To see that this is surjective on objects take an extension
$A \to B$ of $\Lambda$-algebras; from the induced 
${\mcat U}_A \to {\mcat U}_B$ we get ${\mcat U}_A$-algebra
structure on $B \cong {\mcat U}_B(0)$. The functor is evidently faithful.
To see it is full, suppose $B \to C$ is a $\Lambda$-algebra
map between extensions $A \to B$ and $A \to C$ coming from
${\mcat U}_A$-algebras $B$ and $C$; we get ${\mcat U}_B \to {\mcat U}_C$
under ${\mcat U}_A$ and so have a corresponding ${\mcat U}_A$-algebra
map. Since ${\rm Alg}({\mcat U}_A) \to{\rm Alg}({\Lambda}_A)$ 
is an equivalence Proposition \ref{algext} shows that 
${\Lambda}_A \to {\mcat U}_A$ is an isomorphism of algebraic theories.
So in particular there is a canonical $\lambda$-theory
structure on ${\Lambda}_A$.
Putting all that together with Propositions \ref{eta} and \ref{eps} gives
the following.
\begin{theorem}\label{fundthm} 
{\rm (Fundamental Theorem of the $\lambda$-Calculus)}
There is an adjoint equivalence ${\mcat L} \mapsto {\mcat L}(0)$,
$A \to \Lambda_A$ between $\lambda$-theories and
$\Lambda$-algebras: for $\mcat L$ a $\lambda$-theory and 
$A$ a $\Lambda$-algebra,
there are natural isomorphisms ${\mcat L} \cong \Lambda_{{\mcat L}(0)}$ and
$\Lambda_A(0) \cong A$. In particular
each $\lambda$-theory $\mcat L$ is isomorphic 
to the theory of extensions of its initial algebra ${\mcat L}(0)$. 
\end{theorem}
There are many straightforward consequences of the Fundamental Theorem.
Note that for every $\lambda$-theory
$\mcat L$, the ${\mcat L}(n)$ are not just
the free extensions of ${\mcat L}(0)$ as
$\mcat L$-algebras but also as $\Lambda$-algebras.
Much the same thought is expressed in the following.
\begin{proposition} Suppose that $A$ is a $\Lambda$-algebra.
Then there is a canonical $\Lambda$-algebra stucture on the
retracts $A(n) = \{ a \in A | \app_n a = a \}$ 
making them the free $\Lambda$-algebra
extending $A$ by $n$ indeterminates.
\end{proposition}
{\bf Proof} Immediate given the remark 
at the end of Section \ref{fnspace}.\\[0.3em]
This last result is folklore mentioned in
passing in \cite{Freyd89}
and spelled out in \cite{Sel02}. It arises very naturally 
in the approach which I have laid out, and itself suggests
the development of connections with combinatory logic. 

\subsection{An alternative approach}

I can readily imagine that syntactically minded 
colleagues will not be 
comfortable with my approach to the Fundamental
Theorem, so I briefly sketch an alternative.
Let us try to construct $\Lambda_A$
syntactically. Given a $\lambda$-algebra $A$, take 
an extension of the syntax of the $\lambda$-calculus
with constants from $A$. Let $\Lambda_A(n)$ be the terms
with $n$ variables factored out by the equality generated
by $\beta$-equality in the $\lambda$-calculus and by
the equalities given by the actions $\Lambda(m) \times A^m \to A$.
Extending the argument for the initial
$\lambda$-theory $\Lambda$, one can show that the 
resulting $\Lambda_A$ is a $\lambda$-theory. Functoriality of the operation
$A \mapsto \Lambda_A$ is straightforward, but after that
things get delicate as in this approach we do not know that
$\Lambda_A$ is the algebraic theory of extensions of the $\Lambda$-algebra $A$.
Even the isomorphism $A \cong \Lambda_A(0)$ is not obvious:
how do we know that the syntactic theory does not produce
a proper quotient of $A$? There are direct arguments for that
but probably it is easiest to use Proposition \ref{abstraction}
to identify the syntactic $\Lambda_A$ with the theory of extensions.

Next what happens if we start with a $\lambda$-theory $\mcat L$
and form $\Lambda_{{\mcat L}(0)}$? We now know $\Lambda_A$ as
the theory of extensions and so we get a factorization
$\Lambda \to \Lambda_{{\mcat L}(0)} \to {\mcat L}$ of algebraic
theories. By Proposition \ref{theorymap} we deduce that 
$\Lambda_{{\mcat L}(0)} \to {\mcat L}$ is a map of 
$\lambda$-theories. We want to show that this is an
isomorphism but we cannot exploit 
Proposition \ref{algext} as we have no handle on
$\mcat L$-algebras. It seems best to use Proposition \ref{abstraction}
again and argue directly that any $\lambda$-theory $\mcat L$ is the theory
of extensions of the $\Lambda$-algebra ${\mcat L}(0)$. That can be done
but subtleties which appear in Section \ref{fundsec} cannot be
avoided. Overall the syntactic approach is not as straightfoward
as it seems.

\section{Conclusions and Vistas}

I want to stress that
understanding $\lambda$-theories comes before 
understanding $\Lambda$-algebras or any other equivalent notion.
Though the case for the theory approach is foundational, I hope that 
this paper
will encourage new research.
For example customary questions about the syntactic
theories represented by concrete 
cartesian closed categories seem narrow: one should ask 
about the $\lambda$-theories represented. What can one say about
them? How for example to compare them within and between categories?
The Taylor Fibration of Section
\ref{taylor} exhibits for every $\lambda$-theory $\mcat L$
a large family of interpretations corresponding to the
slices of $P{\mcat L}$. But there are many more endomorphism
$\lambda$-theories of reflexive objects in
$P{\mcat L}$. Can one characterize them?

Another particularly interesting set of questions concerns
cartesian closed categories arising from models of the differential lambda
calculus. A form of the Approximation
Theorem holds automatically and it follows that the quotients of
the syntactic theories obtained are very restricted. 
These issues are intimately tied up with B\"{o}hm trees.
For background see \cite{ER03}, \cite{ER06} and \cite{Car11}.
However it appears that notwithstanding the 
restrictions, there is still a wide variety of interpretations.
Is this impression true? How can one make sense of it? 
What general tools are there for telling
differences in such cases?

I close with some remarks about the potential
wider significance of the techniques
discovered by Corrado B\"{o}hm and presented in his
seminal paper \cite{Bohm68}. This is surely a cornerstone of our understanding
of the $\lambda$-calculus and should come early in any account
of fundamental ideas. But how well do we understand what is involved?
Prima facie the techniques are syntactic. The original applications are to 
quotients of the initial $\beta\eta$ theory $\Lambda_{\eta}$ and concern
what are usually thought of as the limits of consistency: what can you or 
can you not do before such a quotient becomes the terminal 
(trivial) theory? There is B\"{o}hm's original point that
identifying distinct $\beta\eta$ normal forms collapses a theory
and the closely related result of \cite{Hyl76} that there is a unique
maximal non-trivial quotient of the theory extending  $\Lambda_{\eta}$
by setting all unsolvables equal. But I think there is a broader set
of semantic principles at stake. For example when analysing models such as 
Scott's $P\omega$ in which $\eta$ does not hold one exploits
non-syntactic variants of B\"{o}hm's combinators. (This is well
explained for Plotkin's
$T^{\omega}$ model
in section 3 of \cite{BL80}.) 
There must surely be more to be said on a broader 
canvas.

In honouring B\"{o}hm I want to stress the future.
There is still much to discover about the 
pure $\lambda$-calculus. Here I have presented it in modern
dress with all the trappings of current categorical
research. My intention is to convey the clear message 
that the $\lambda$-calculus is a remarkable 
and exceptionally elegant area of abstract mathematics.

\end{document}